\documentclass[11pt]{amsart}
\usepackage{amsmath,amsthm,epsfig,amssymb}
\usepackage{xypic}
\input xy
\title{Ehrhart's polynomial for equilateral triangles in $\mathbb Z^3$}
\author{Eugen J. Ionascu}
\curraddr{Department of Mathematics\\ Columbus State University\\4225 University Avenue\\
Columbus, GA 31907\\
Honorific Member of the Romanian Institute of Mathematics ``Simion
Stoilow" } \email{ionascu@columbusstate.edu;} 
\date{July $4^{th}$, 2011}
\begin{document}
\def\sms{\small\scshape}
\baselineskip18pt
\newtheorem{theorem}{\hspace{\parindent}
T{\scriptsize HEOREM}}[section]
\newtheorem{proposition}[theorem]
{\hspace{\parindent }P{\scriptsize ROPOSITION}}
\newtheorem{corollary}[theorem]
{\hspace{\parindent }C{\scriptsize OROLLARY}}
\newtheorem{lemma}[theorem]
{\hspace{\parindent }L{\scriptsize EMMA}}
\newtheorem{definition}[theorem]
{\hspace{\parindent }D{\scriptsize EFINITION}}
\newtheorem{problem}[theorem]
{\hspace{\parindent }P{\scriptsize ROBLEM}}
\newtheorem{conjecture}[theorem]
{\hspace{\parindent }C{\scriptsize ONJECTURE}}
\newtheorem{example}[theorem]
{\hspace{\parindent }E{\scriptsize XAMPLE}}
\newtheorem{remark}[theorem]
{\hspace{\parindent }R{\scriptsize EMARK}}
\renewcommand{\thetheorem}{\arabic{section}.\arabic{theorem}}
\renewcommand{\theenumi}{(\roman{enumi})}
\renewcommand{\labelenumi}{\theenumi}
\newcommand{\Q}{{\mathbb Q}}
\newcommand{\Z}{{\mathbb Z}}
\newcommand{\N}{{\mathbb N}}
\newcommand{\C}{{\mathbb C}}
\newcommand{\R}{{\mathbb R}}
\newcommand{\F}{{\mathbb F}}
\newcommand{\K}{{\mathbb K}}
\newcommand{\D}{{\mathbb D}}
\def\phi{\varphi}
\def\ra{\rightarrow}
\def\sd{\bigtriangledown}
\def\ac{\mathaccent94}
\def\wi{\sim}
\def\wt{\widetilde}
\def\bb#1{{\Bbb#1}}
\def\bs{\backslash}
\def\cal{\mathcal}
\def\ca#1{{\cal#1}}
\def\Bbb#1{\bf#1}
\def\blacksquare{{\ \vrule height7pt width7pt depth0pt}}
\def\bsq{\blacksquare}
\def\proof{\hspace{\parindent}{P{\scriptsize ROOF}}}
\def\pofthe{P{\scriptsize ROOF OF}
T{\scriptsize HEOREM}\  }
\def\pofle{\hspace{\parindent}P{\scriptsize ROOF OF}
L{\scriptsize EMMA}\  }
\def\pofcor{\hspace{\parindent}P{\scriptsize ROOF OF}
C{\scriptsize ROLLARY}\  }
\def\pofpro{\hspace{\parindent}P{\scriptsize ROOF OF}
P{\scriptsize ROPOSITION}\  }
\def\n{\noindent}
\def\wh{\widehat}
\def\eproof{$\hfill\bsq$\par}
\def\ds{\displaystyle}
\def\du{\overset{\text {\bf .}}{\cup}}
\def\Du{\overset{\text {\bf .}}{\bigcup}}
\def\b{$\blacklozenge$}

\def\eqtr{{\cal E}{\cal T}(\Z) }
\def\eproofi{\bsq}

\begin{abstract} In this paper we calculate the Ehrhart's
polynomial associated with a 2-dimensional regular polytope
(i.e. equilateral triangles) in $\mathbb Z^3$. The polynomial takes a relatively simple form in terms of the
coordinates of the vertices of the polytope and it depends heavily on the value $d$ and its divisors, where
$d=\sqrt{\frac{a^2+b^2+c^2}{3}}$ and $(a,b,c)$ ($\gcd(a,b,c)=1$) is a vector with integer coordinates normal to the plane containing the triangle.
\end{abstract} \maketitle
\section{INTRODUCTION}
A description of all equilateral triangles with vertices in
$\mathbb Z^3$ appeared first in \cite{eji} (with the proof of the
full general case in \cite{rceji}). An updated version of the same
results but with a shorter analysis was included in
\cite{ejiobando}.

In the 1960's, $\rm Eug\grave{e}ne$ Ehrhart (\cite{ee},\cite{ee2}) proved that given a $\mathbb Z^n$ lattice $d$-dimensional polytope in $\mathbb
R^n$ ($1\le d\le n$), denoted here generically by $\cal P$,
there exists a polynomial  $L({\cal P}, t)\in \mathbb
Z[t]$, associated with $\cal P$, satisfying

\begin{equation} \label{ehrhartpolynomial}
L({\cal P},t)=the \ cardinality\  of\ \{t\cal P\}\cap \mathbb Z^n, \ t\in \mathbb N.
\end{equation}

\n Lately there has been a lot of literature on various topics generated by
the Ehrhart polynomial and we cite just a few papers and books:
\cite{ChristosAAthanasiadis}, \cite{abarvinok}, \cite{beckAndRobins2007textbook}, \cite{diazbeckrobi}, \cite{diazrobins1997},  \cite{BJBraun}, \cite{BJBraun2},    \cite{BJBraun3},\cite{BrionVergne}, \cite{BChen}, \cite{FuLiu}, \cite{StevenSam}, and \cite{AlanStapledon}.

Every equilateral triangle in $\mathbb Z^3$
after a translation by a vector with integer coordinates can be
assumed to have the origin as one of its vertices. Then one can
show that the triangle's other vertices are contained in a lattice
of points of the form (see \cite{eji})

\begin{equation}\label{planelattice}
{\cal P}_{a,b,c}:=\{(x,y,z)\in \mathbb {Z}^3|\ \
ax+by+cz=0,\ \ a^2+b^2+c^2=3d^2,\ \ gcd(a,b,c)=1, \
a,b,c,d\in \mathbb Z\}.
\end{equation}

\begin{center}\label{figure0}
\epsfig{file=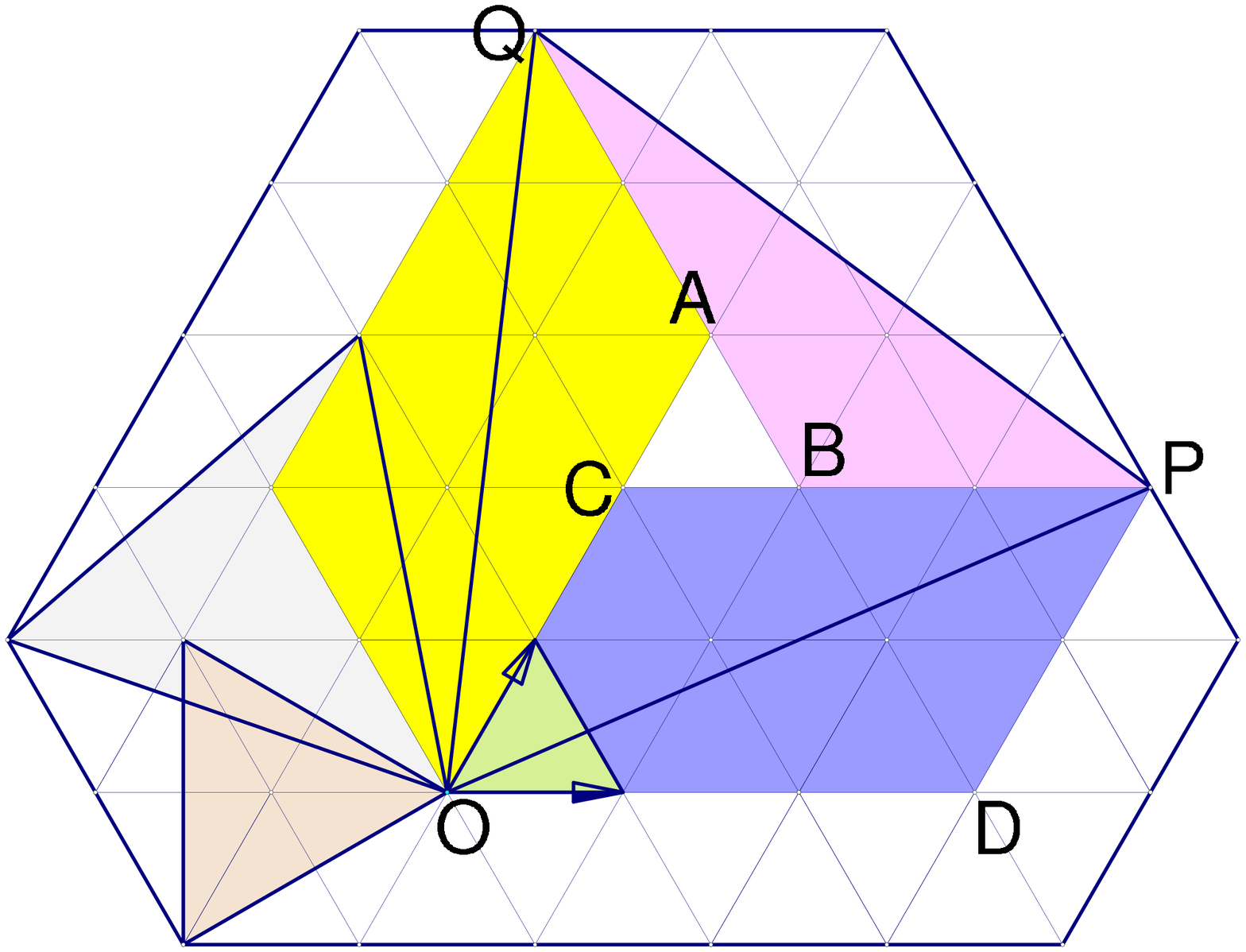,height=2in,width=2.3in}
\vspace{.1in}\\
{\small Figure 1: The lattice ${\cal P^{eq}}_{a,b,c}$ } \vspace{.1in}
\end{center}

This lattice is in general much richer than the sub-lattice,
${\cal P^{eq}}_{a,b,c}$, of all points which are vertices of
equilateral triangles with one of the vertices the origin (see
Figure~\ref{figure0}). A more precise description of this
sub-lattice is given by the following result essentially
contained in \cite{ejiobando}.

\begin{theorem}\label{oldnewparametrization}
The sub-lattice ${\cal P^{eq}}_{a,b,c}$ is generated by two
vectors $\overrightarrow{\zeta}$ and $\overrightarrow{\eta}$ in
the following sense: ${{\cal T}^{m,n}_{a,b,c}}:=\triangle OPQ$ with
$P$, $Q$ in ${\cal P}_{a,b,c}$, is equilateral if and only if for
some integers $m$, $n$

\begin{equation}\label{vectorid}
\overrightarrow{OP}=m\overrightarrow{\zeta}-n\overrightarrow{\eta},\
\
\overrightarrow{OQ}=n\overrightarrow{\zeta}+(m-n)\overrightarrow{\eta},
\ { \rm with} \ \overrightarrow{\zeta}=(\zeta_1,\zeta_2,\zeta_3),
\overrightarrow{\varsigma}=(\varsigma_1,\varsigma_2,\varsigma_3),
\overrightarrow{\eta}=\frac{\overrightarrow{\zeta}+\overrightarrow{\varsigma}}{2},
\end{equation}

\begin{equation}\label{paramtwo}
\begin{array}{l}
\begin{cases}
\ds \zeta_1=-\frac{rac+dbs}{q} \\ \\
\ds \zeta_2=\frac{das-bcr}{q}\\ \\
\ds \zeta_3=r
\end{cases}
\ \ ,\ \ \ \begin{cases}
\ds \varsigma_1=\frac{3dbr-acs}{q} \\ \\
\ds \varsigma_2=-\frac{3dar+bcs}{q}\\ \\
\ds \varsigma_3=s
\end{cases}
\end{array}\ \ ,
\end{equation}
where $q=a^2+b^2$ and $r$, $s$ can be chosen so that all six
numbers in (\ref{paramtwo}) are integers. The sides-lengths of
$\triangle OPQ$ are equal to $d\sqrt{2(m^2-mn+n^2)}$. Moreover,
$r$ and $s$ can be selected in such a way that the following
properties are also true: \par

\n (i) $r$ and $s$ satisfy  $2q=s^2+3r^2$ and as a result
$2(b^2+c^2)=\varsigma_1^2+3\zeta_1^2$ and
$2(a^2+c^2)=\varsigma_2^2+3\zeta_2^2$\par

\n (ii) $r=r'\omega \chi$, $s=s'\omega\chi $ where
$\omega=gcd(a,b)$, $gcd(r',s')=1$ and $\chi$ is the product of the
prime factors of the form $6k-1$ of $gcd(d,q)$

\n (iii) $\chi$ divides $c$

\n (iv) $|\overrightarrow{\zeta}|=d\sqrt{2}$,
$|\overrightarrow{\varsigma}|=d\sqrt{6}$, and
$\overrightarrow{\zeta}\cdot\overrightarrow{\varsigma} =0$.

\n (v) $s+i\sqrt{3}r=gcd(A-i\sqrt{3}B,2q)$, in the ring $\mathbb Z[i\sqrt{3}]$, where
$A=ac$ and $B=bd$

\n (vi) This construction is essentially independent of the choice of $a$ and
$b$. In other words, similar statements can be formulated for $a$
and $c$, or $b$ and $c$ instead of $a$ and $b$.
\end{theorem}

\n {\bf Remark:}  From a computational point of view, it is better to use the smaller two of the numbers $a$, $b$ and $c$.

An example here may be illuminating. If $d=15$, we observe that we
can take $a=1$, $b=7$ and $c=25$ ($3d^2=a^2+b^2+c^2$). Then,
$r=-5$ and $s=5$ give $\overrightarrow{\zeta}=(13,16,-5)$ and
$\overrightarrow{\eta}=(21,-3,0)$. Also,  properties (i), (ii),
and (iii) in Theorem~\ref{oldnewparametrization} are satisfied.

In this paper we are interested in finding the Ehrhard polynomial for this type of triangle in $\mathbb Z^3$.
Our formulae seem to be simpler than various other formulae obtained for other polytopes. This may be caused by the symmetries
and the complications that are already built into our polytopes (see the precise formulation below (\ref{theproblem})).

From the general theory (see \cite{beckAndRobins2007textbook} for a good account) we know that in this case

$$L({\Delta },t)=c_0t^2+\frac{c_1}{2}t+c_2, \ t\in \mathbb N,$$

\n  where $c_2=1$ since we are dealing with a (convex) polytope and $c_0$ is the area of the triangle $\Delta$ normalized to
the area of a fundamental domain of the sub-lattice ${\cal P}_{a,b,c}$. We will see that $c_0$ is really easy to compute but
$c_1$, the number of points of the sub-lattice ${\cal P}_{a,b,c}$ on the sides of the triangle, creates some difficulties.
Let us point out that this polynomial is not a complete invariant, in the sense that we may have the same
Ehrhart polynomial for two ``different" triangles. By different triangles, we understand that one triangle cannot be obtained
from the other by the usual transformations which leave the lattice
  $\mathbb Z^3$ invariant. If we take $\Delta_1:=\{(0,0,0),(13,-8,3),(0,-11,11)\}$ and $\Delta_2:=\{(0,0,0),(4,15,-1),(15,4,-1)\}$,
then

$$ L({\Delta }_1,t)=L({\Delta }_2,t)=\frac{11}{2}t^2+\frac{13}{2}t+1, \ t\in \mathbb N.$$
The triangles are essentially different since they live in totally different sub-lattices: ${\cal P}_{5,13,13}$ and  ${\cal P}_{1,1,19}$ respectively.

It turns out that what we are calculating reduces to counting the number of integer triples $(x,y,z)$ satisfying

\begin{equation}\label{theproblem}
\begin{cases}
ax+by+cz=0\\ \\
(3dbr-acs)x-(3dar+bcs)y+qsz\ge 0\\ \\
[ac(s-3r)-3db(r+s)x+[3da(r+s)+bc(s-3r)]y+(3r-s)qz\ge 0\\ \\
[db(r+s)+ac(s-r)]x-[da(3r-s)+bc(r+s)]y+q(s-r)z+2qd^2t^2\ge 0
\end{cases}, t\in \mathbb N.
\end{equation}

\vspace{0.2in}

\section{A fundamental domain in ${\cal P}_{a,b,c}$ and the coefficient $c_0$ }
\vspace{0.2in}

In what follows we are going to assume that $gcd(a,b,c)=1$ and use
the notation introduced in Theorem~\ref{oldnewparametrization}.
We notice that from the relations (\ref{vectorid}) and
(\ref{paramtwo}), we obtain

\begin{equation}\label{interiororbit}
\frac{r+s}{2} \overrightarrow{\zeta}-r\overrightarrow{\eta}=d
(-b,a,0), \ \frac{\zeta_1+\varsigma_1}{2}
\overrightarrow{\zeta}-\zeta_1\overrightarrow{\eta}=d (0,-c,b), \
and\ \frac{\zeta_2+\varsigma_2}{2}
\overrightarrow{\zeta}-\zeta_2\overrightarrow{\eta}=d (c,0,-a).
\end{equation}

 \n It is clear that the vectors
$\overrightarrow{u}=\frac{1}{\omega}(-b,a,0)$,
$\overrightarrow{v}=\frac{1}{gcd(a,c)}(-c,0,a)$ and
$\overrightarrow{w}=\frac{1}{gcd(b,c)}(0,-c,b)$ correspond to
points in ${\cal P}_{a,b,c}$. We first show that
${\cal P}_{a,b,c}$ is generated by $\overrightarrow{u}$,
$\overrightarrow{v}$ and $\overrightarrow{w}$ (a point $P$ in
${\cal P}_{a,b,c}$ is identified by its position vector
$\overrightarrow{OP}$ as usual).
\begin{lemma}\label{basisofp} With the above notation, we have $${\cal
P}_{a,b,c}=\{m\overrightarrow{u}+n\overrightarrow{v}+p\overrightarrow{w}|m,n,p\in
\mathbb Z\}. $$
\end{lemma}
\n \proof. If $(x,y,z)\in {\cal P}_{a,b,c}$, we see
that $ax+by+cz=0$. Because $\omega=gcd(a,b)$ and
$gcd(a,b,c)=1$ we need to have $z=\omega z'$ with
$z'\in \mathbb Z$. Also, the existence of integers $k$ and
$l$ such that $ka+\ell b=\omega$ is insured by the fact
$\omega=gcd(a,b)$. This means that we have
$$(x,y,z)-z'[gcd(a,c)k\overrightarrow{v}+gcd(b,c)\ell \overrightarrow{w}]=(\alpha',\beta',0)\in {\cal P}_{a,b,c}.$$
Since $a\alpha'+b\beta'=0$ we see that $(\alpha',\beta',0)=\lambda
\overrightarrow{u}$ for some $\lambda\in \mathbb  Z^3$. This shows
the inclusion $\subset$ in the equality claimed, and the inclusion
$\supset$ is obvious.\eproof

If we look at the proof of the above lemma we see that it is
not necessary to have three vectors to generate ${\cal
P}_{a,b,c}$, but only $\overrightarrow{u}$ and
$\overrightarrow{\tau}:=gcd(a,c)k\overrightarrow{v}+gcd(b,c)\ell
\overrightarrow{w}$, where $ka+\ell b=gcd(a,b)$, are enough. Since there are infinitely many pairs $(k,l)$ satisfying this equality, let us take the solution that minimizes $k$ so that $k>0$.    For computational purposes, we have the following more useful result.

\begin{lemma}\label{basisofp2} With the above definition of $\overrightarrow{\tau}$ we have

(i) ${\cal
P}_{a,b,c}=\{k\overrightarrow{u}+\ell \overrightarrow{\tau}|\ k,\ell\in
\mathbb Z\}$,

(ii) $\overrightarrow{u}=\frac{\widetilde{r}+\widetilde{s}}{2 d} \overrightarrow{\zeta}-\frac{\widetilde{r}}{ d}\overrightarrow{\eta}$, where $r=\omega\widetilde{r}$ and $s=\omega\widetilde{s}$ (Theorem~\ref{oldnewparametrization}) and

\begin{equation}\label{secondvector}
\overrightarrow{\tau}=\frac{\alpha}{d} \overrightarrow{\zeta}+\frac{\beta}{d}\overrightarrow{\eta},\ \ for\ some\  \alpha,\beta\in \mathbb Z.
\end{equation}
\end{lemma}

\n \proof. The first part follows from Lemma~\ref{basisofp} and (ii) is a consequence of the equalities (\ref{interiororbit}).\eproof

\n So, $\overrightarrow{u}$ and $\overrightarrow{\tau}$ form a fundamental domain for
${\cal
P}_{a,b,c}$. The area of the parallelogram formed by these vectors is given by
$|\overrightarrow{u}\times \overrightarrow{\tau}|$.

\begin{lemma}\label{areaoffundamentaldomain} (i) The area of a fundamental domain for ${\cal
P}_{a,b,c}$ is equal to $d\sqrt{3}$.

(ii) The integers $\alpha$ and $\beta$ in (\ref{secondvector}) satisfy the relation

\begin{equation}\label{alphabeta} |\frac{\widetilde{r}+\widetilde{s}}{2}\beta+\widetilde{r}\alpha|=d
\end{equation}
\end{lemma}

\n \proof. (i) We observe that $\overrightarrow{u}\times \overrightarrow{v}=\frac{a}{\omega \gcd(a,c)}(a\overrightarrow{i}+b\overrightarrow{j}+c\overrightarrow{k})$
and similarly $\overrightarrow{u}\times \overrightarrow{w}=\frac{b}{\omega \gcd(b,c)}(a\overrightarrow{i}+b\overrightarrow{j}+c\overrightarrow{k})$.
Hence the area of the parallelogram determined by $\overrightarrow{u}$ and $\overrightarrow{\tau}$ is equal to

$$|\overrightarrow{u}\times \overrightarrow{\tau}|=|\frac{1}{\omega}(ak+bl)(a\overrightarrow{i}+b\overrightarrow{j}+c\overrightarrow{k})|=\sqrt{a^2+b^2+c^2}=d\sqrt{3}.$$

The second part follows from (i) and the equality (\ref{secondvector}), if we take into account that, $\overrightarrow{\zeta}\times \overrightarrow{\zeta}=\overrightarrow{\eta}\times \overrightarrow{\eta}=0$,
$\overrightarrow{\zeta}\times \overrightarrow{\eta}=-\overrightarrow{\eta}\times \overrightarrow{\zeta}$, and
$|\overrightarrow{\zeta}\times \overrightarrow{\eta}|=d^2\sqrt{3}$:

$$d\sqrt{3}=|\overrightarrow{u}\times \overrightarrow{\tau}|=|\frac{\frac{\widetilde{r}+\widetilde{s}}{2}\beta+\widetilde{r}\alpha}{d^2}\overrightarrow{\zeta}\times \overrightarrow{\eta}|=\frac{|\frac{\widetilde{r}+\widetilde{s}}{2}\beta+\widetilde{r}\alpha|}{d^2}d^2\sqrt{3}\Rightarrow (\ref{alphabeta}).$$

\begin{proposition}\label{coefficientc0}
The coefficient $c_0$ in the Ehrhart polynomial associated with an equilateral triangle
${{\cal T}^{m,n}_{a,b,c}}$ described by Theorem~\ref{oldnewparametrization} is given by
$$c_0=\frac{d(m^2-mn+n^2)}{2}.$$
\end{proposition}

\n \proof. By the general theory of the Ehrhart polynomial, $c_0$ is equal to the area of the triangle normalized by the area of a fundamental domain
of the lattice ${\cal P}_{a,b,c}$. Since the area of the triangle ${{\cal T}^{m,n}_{a,b,c}}$ is equal to $\frac{2d^2(m^2-mn+n^2)\sqrt{3}}{4}$, using Lemma~\ref{areaoffundamentaldomain},
we obtain $c_0=\frac{2d^2(m^2-mn+n^2)\sqrt{3}}{4d\sqrt{3}}=\frac{d(m^2-mn+n^2)}{2}$. \eproof

\section{The coefficient $c_1$ }

\begin{center}\label{figure2}
$\epsfig{file=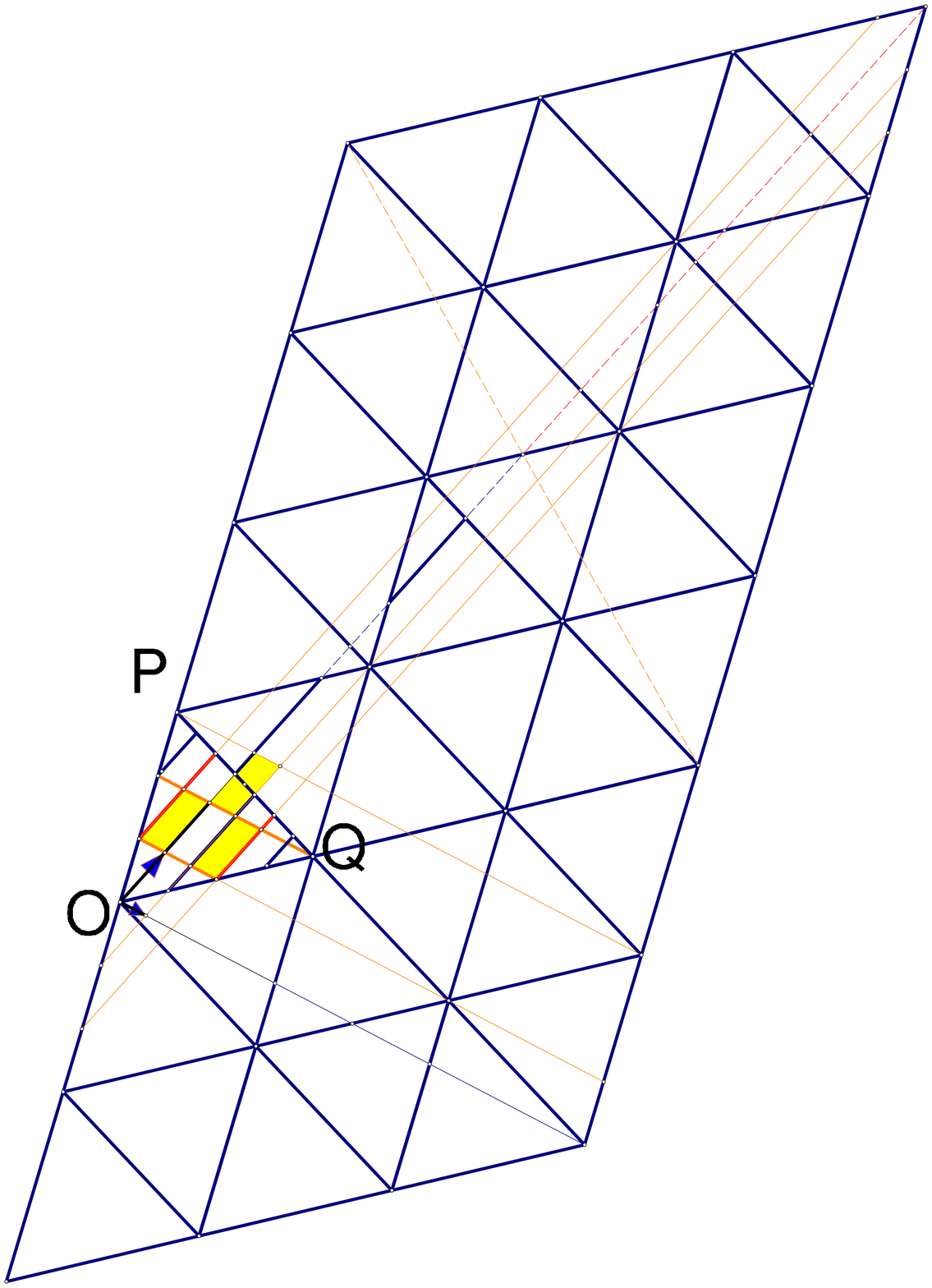,height=2.5in,width=2.7in}$
\vspace{.1in}\\
{\small Figure 2:
$\overrightarrow{u}=\frac{\widetilde{r}+\widetilde{s}}{2 d}\overrightarrow{\zeta}-\frac{\widetilde{r}}{d}\overrightarrow{\eta},\ \ \overrightarrow{\tau}=\frac{\alpha}{d} \overrightarrow{\zeta}+\frac{\beta}{d}\overrightarrow{\eta}$}
\end{center}

 We consider the case $m=1$ and $n=0$. From Lemma~\ref{basisofp}, we see (Figure~2) that the problem at hand is to count the points of the lattice, that are determined by  $\overrightarrow{u}$ and $\overrightarrow{\tau}$, on the sides of the triangle  ${{\cal T}^{1,0}_{a,b,c}}$.  Clearly we have $3$ points, the vertices of the triangle,
and in most cases these are all of the lattice points on the sides ($c_1=3$).
The points on the sides can be (see Figure~2) on the side $\overline{OP}$, $\overline{PQ}$ or $\overline{OQ}$.
We note that the points on the side $\overline{OP}$, are characterized by the existence of integers $k$, $\ell$, $\lambda$ and $\mu$ such that

\begin{equation}\label{sideOP} k\overrightarrow{u}+\ell \overrightarrow{\tau}=(k\frac{\widetilde{r}+\widetilde{s}}{2 d}+\ell \frac{\alpha}{d})\overrightarrow{\zeta}+(-k\frac{\widetilde{r}}{d}+\ell\frac{\beta}{d})\overrightarrow{\eta}=(k\frac{\widetilde{r}+\widetilde{s}}{2d}+\ell \frac{\alpha}{d}) \overrightarrow{\zeta}+\mu\overrightarrow{\eta}.
\end{equation}

Similarly, for points on the side $\overline{OQ}$, the equality above changes to

\begin{equation}\label{sideOQ} k\overrightarrow{u}+\ell \overrightarrow{\tau}=(k\frac{\widetilde{r}+\widetilde{s}}{2 d}+\ell \frac{\alpha}{d})\overrightarrow{\zeta}+(-k\frac{\widetilde{r}}{d}+\ell\frac{\beta}{d})\overrightarrow{\eta}=\lambda \overrightarrow{\zeta}+ (-k\frac{\widetilde{r}}{d}+\ell\frac{\beta}{d})\overrightarrow{\eta}.\end{equation}

For the side $\overline{PQ}$, we have the characterization

\begin{equation}\label{sidePQ} k\frac{\widetilde{s}-\widetilde{r}}{2d}+\ell\frac{\alpha+\beta}{d}\in \mathbb Z.
\end{equation}

\begin{lemma}\label{pointsonTheSides}
(i) The coefficient $(-k\frac{\widetilde{r}}{d}+\ell\frac{\beta}{d})$ is an integer if $\ell =\mu\frac{\widetilde{r}+\widetilde{s}}{2}+ \lambda \widetilde{r}$ and $k=\lambda \beta-\mu\alpha$ for every
$\lambda, \mu \in \mathbb Z^3$.

(ii)  For the given values of $\ell$ and $k$ above the other coefficient in (\ref{sideOP}) becomes

$$k\frac{\widetilde{r}+\widetilde{s}}{2 d}+\ell \frac{\alpha}{d}=\pm  \lambda\in \mathbb Z.$$
(iii) If $\gcd(\widetilde{r},\beta)=1$, then there are no  points of the lattice in the interior  of $\overline{OP}$.
\end{lemma}

\proof. Using (\ref{alphabeta}), we have
$$(-k\frac{\widetilde{r}}{d}+\ell\frac{\beta}{d})=(-\lambda \beta+\mu\alpha)\frac{\widetilde{r}}{d}+ (\mu\frac{\widetilde{r}+\widetilde{s}}{2}+ \lambda \widetilde{r})\frac{\beta}{ d}=\pm \mu.$$
Similarly, one checks that (ii) is true.

\n We observe that (\ref{alphabeta}) implies in particular that $\gcd(\widetilde{r},\beta)$ divides $d$. Therefore, if  $\gcd(\widetilde{r},\beta)=1$, the values in (i) for $k$ and $\ell$ give all of the solutions of the Diophantine equation
$(-k\frac{\widetilde{r}}{d}+\ell\frac{\beta}{d})=\mu\in \mathbb Z$. (We refer the reader to a basic text on linear Diophantine equations such as \cite{mnz} or \cite{r}.) Part (ii) implies that there are no lattice points on the side $\overline{OP}$ other than the endpoints. \eproof

Clearly, a similar lemma is true for the sides $\overline{OQ}$ and $\overline{PQ}$.

So, we have proved the following proposition.

\begin{proposition}\label{firstresult} If  $\gcd(\widetilde{r},\beta)=\gcd(\frac{\widetilde{r}+\widetilde{s}}{2},\alpha)=\gcd(\frac{\widetilde{r}-\widetilde{s}}{2},\frac{\alpha+\beta}{2})=1$ then there are no lattice points on the sides of ${{\cal T}^{1,0}_{a,b,c}}$ other than the vertices. The Ehrhart polynomial in this case is
$$L({\cal T}^{1,0}_{a,b,c},t)=\frac{dt^2+3t}{2}+1, \ t\in \mathbb N.$$
\end{proposition}

\n It is natural to ask if all sides of the equilateral triangle ${\cal T}^{1,0}_{a,b,c}$ may simultaneously contain lattice points in their interiors and how many can there be?
The next lemma answers these questions.

\begin{center}\label{figure3}
$$\underset{\small Figure\ 3:\
All\ sides\ may\ contain\ lattice \ points}{\epsfig{file=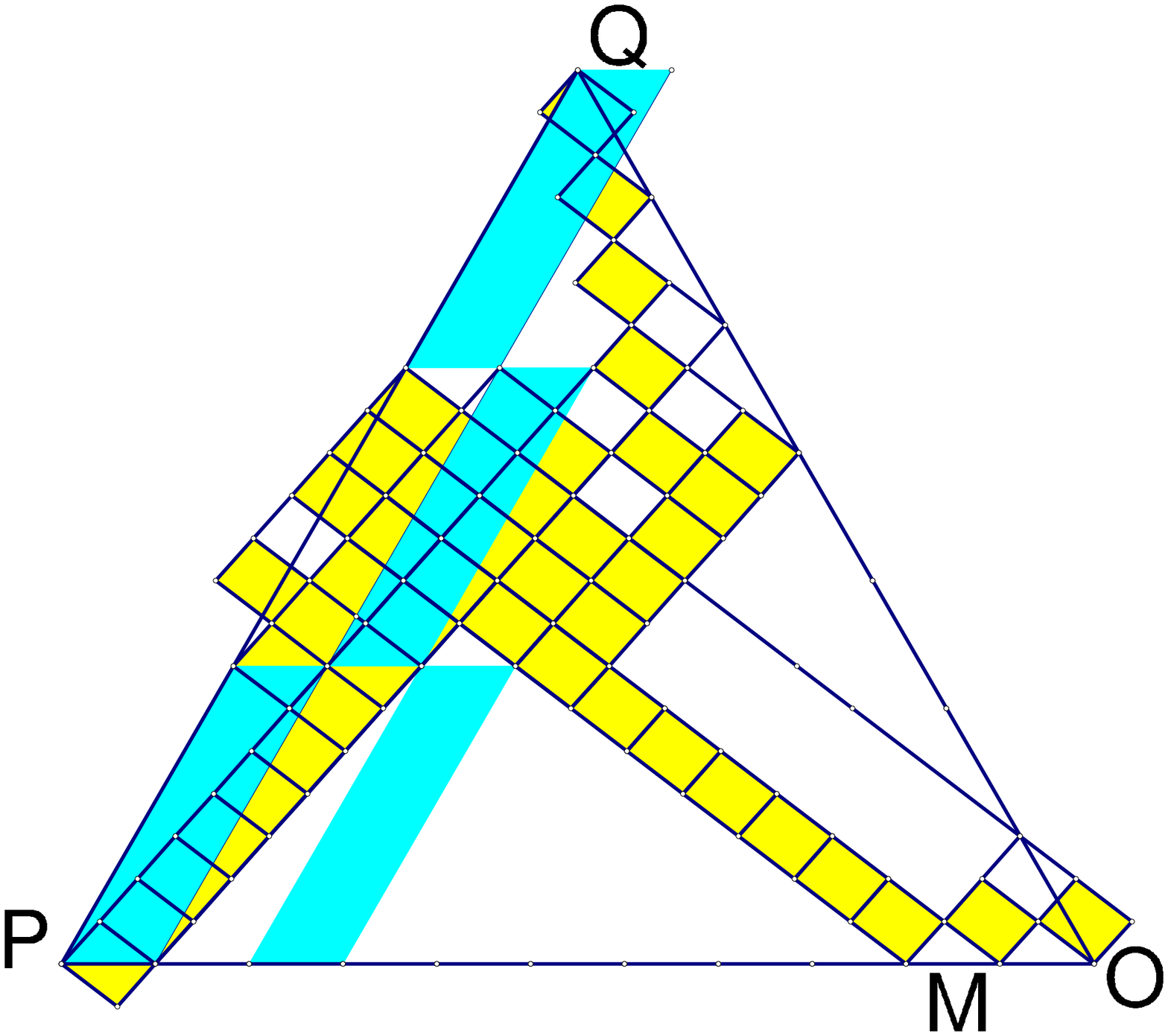,height=2in,width=2in}}$$
\end{center}

\begin{lemma}\label{noallsides} (i) If a side contains lattice points in its interior, the intersection of that side with ${\cal P}_{a,b,c}$ is the set of points which divide that side into $d'$ equal parts, where $d'$ divides $d$.\par

(ii) It is possible that all of the sides of a minimal triangle ${\cal T}^{1,0}_{a,b,c}$  have lattice points of ${\cal P}_{a,b,c}$ in their interiors. If $d_1$ and $d_2$ are the corresponding numbers as in part (i), for two sides, then $\gcd(d_1,d_2)=1$.
\end{lemma}

\proof. (i) We refer to Figure~3 in this proof. Without loss of generality, let us assume that the closest point to $O$ in ${\cal P}_{a,b,c}$ on the side $\overline{OP}$, is $M$.
Then $\overrightarrow{OM}$ has integer coordinates. Then so does $k\overrightarrow{OM}$ for all $k\in \mathbb Z$. Therefore, for some $k$ we must have $k\overrightarrow{OM}=\overrightarrow{OP}$, otherwise we can construct some point in   ${\cal P}_{a,b,c}$, on the side  $\overline{OP}$, closer to $O$ than $M$, using the idea of division with non-zero remainder. Therefore, we have shown the first part of (i).
For the second part we observe that the Diophantine equation $\frac{2d^2}{k^2}=|\overrightarrow{OM}|^2\in \mathbb Z$ is possible only if $k$ divides $d$. This shows (i).

For (ii), we include here the first example we found of such a triangle, $$\Delta_3=\{(0,0,0),(220, 539, -539), (747, 12, -267)\},$$
 \n corresponding to $d=561=3(11)(17)$, $a=245$, $b=613$ and $c=713$. There are respectively 2, 10 and 16 lattice points interior to its sides.  The associated Ehrhart polynomial is

$$L(\Delta_3,t)=\frac{561t^2+31t}{2}+1, \  t\in \mathbb N.$$
For the second part of (ii), let us assume by way of contradiction, that $\gcd(d_1,d_2)=\delta>1$. Then there exist integers $d_1'$ and $d_2'$ such that $d_1=\delta d_1'$ and $d_2=\delta d_2'$. From this we get that $1>\frac{1}{\delta}=\frac{d_1'}{d_1}=\frac{d_2'}{d_2}$. Without loss of generality we can assume that the side $\overline{OP}$ is divided into $d_1$ equal parts and the side $\overline{OQ}$ into $d_2$ equal parts. Then by taking the division points on the sides corresponding to $d_1'$ and $d_2'$ respectively, but different from $O$, we can  construct an equilateral triangle in ${\cal P}_{a,b,c}$, with side lengths
$\frac{d\sqrt{2}}{\delta}$ which are strictly smaller than those of ${\cal T}^{1,0}_{a,b,c}$. Theorem~\ref{oldnewparametrization} shows that this is not possible. The contradiction shows that
$\gcd(d_1,d_2)=1$. \eproof

\n Let us now look at the case $\gcd(\widetilde{r},\beta)=\nu>1$.  Then Lemma~\ref{pointsonTheSides} easily changes  into the following.

\begin{lemma}\label{pointsonTheSidesreally}
(i) If $\gcd(\widetilde{r},\beta)=\nu$, the coefficient $(-k\frac{\widetilde{r}}{d}+\ell\frac{\beta}{d})$ is an integer if $\ell =\mu\frac{\widetilde{r}+\widetilde{s}}{2}+ \lambda \frac{\widetilde{r}}{\nu}$ and $k=\lambda\frac{\beta}{\nu}-\mu\alpha$ for every
$\lambda, \mu \in \mathbb Z^3$.

(ii)  For the given values of $\ell$ and $k$ above the other coefficient in  (\ref{sideOP})  becomes

$$k\frac{\widetilde{r}+\widetilde{s}}{2 d}+\ell \frac{\alpha}{d}=\pm  \frac{\lambda}{\nu} \in \mathbb Z.$$
(iii) If $\gcd(\widetilde{r},\beta)=\nu$, then there are $\nu-1$ points of the lattice in the interior of the side $\overline{OP}$.
\end{lemma}

\proof. The calculation we have done in showing Lemma~\ref{pointsonTheSides} are valid here with the substitution $\frac{\lambda}{\nu}$ instead of $\lambda$. Part (iii) follows from parts (i) and (ii).\eproof

Finally we can put all these things together into the following generalization of Proposition~\ref{firstresult}:

\begin{theorem}\label{coefficientc0}  The Ehrhart polynomial for ${{\cal T}^{1,0}_{a,b,c}}$ (given by Theorem~\ref{oldnewparametrization}) is
$$L({\cal T}^{1,0}_{a,b,c},t)=\frac{dt^2+c_1t}{2}+1, \ t\in \mathbb N,\ where $$
$$3d^2=a^2+b^2+c^2,\ \  c_1=gcd(\widetilde{r},\beta)+\gcd(\frac{\widetilde{r}+\widetilde{s}}{2},\alpha)+\gcd(\frac{\widetilde{r}-\widetilde{s}}{2},\alpha+\beta),$$

\n and ($\alpha$, $\beta$) is a particular solution of the Diophantine equation $\frac{\widetilde{r}+\widetilde{s}}{2}\beta+\widetilde{r}\alpha=d$.
\end{theorem}

We observe that this theorem does not depend on the particular solution $(\alpha, \beta)$ of the Diophantine equation $\frac{\widetilde{r}+\widetilde{s}}{2}\beta+\widetilde{r}\alpha=d$ (or equivalently $\frac{-\widetilde{r}+\widetilde{s}}{2}\beta+\widetilde{r}(\alpha+\beta)=d$).

\section{The general case}
We fix two integers $m$ and $n$ with $\gcd(m,n)=1$ in Theorem~\ref{oldnewparametrization}.
In this case, from (\ref{vectorid}), we have
$\overrightarrow{OP}=m\overrightarrow{\zeta}-n\overrightarrow{\eta}$ and
$\overrightarrow{OQ}=n\overrightarrow{\zeta}+(m-n)\overrightarrow{\eta}$.
We can solve for $\overrightarrow{\zeta}$ and $\overrightarrow{\eta}$:

$$\overrightarrow{\zeta}=\frac{(m-n)\overrightarrow{OP}+n\overrightarrow{OQ}}{m^2-mn+n^2}\ and \ \overrightarrow{\eta}=\frac{m\overrightarrow{OQ}-n\overrightarrow{OP}}{m^2-mn+n^2}.$$

Then, the equations (\ref{sideOP}) and   (\ref{sideOQ}) change into

\begin{equation}
\begin{array}{l}\label{sideOPnew} k\overrightarrow{u}+\ell \overrightarrow{\tau}=(k\frac{\widetilde{r}+\widetilde{s}}{2 d}+\ell \frac{\alpha}{d})\overrightarrow{\zeta}+(-k\frac{\widetilde{r}}{d}+\ell\frac{\beta}{d})\overrightarrow{\eta}= \\ \\
\displaystyle \left(k\frac{m(\widetilde{r}+\widetilde{s})+n(\widetilde{r}-\widetilde{s})}{2d(m^2-mn+n^2)}+\ell \frac{m\alpha-n(\alpha+\beta)}{d(m^2-mn+n^2)}\right)\overrightarrow{OP}+\\ \\ \displaystyle \left(k\frac{n(\widetilde{r}+\widetilde{s})-2m\widetilde{r}}{2d(m^2-mn+n^2)}+\ell \frac{m\beta+n\alpha}{d(m^2-mn+n^2)}\right)\overrightarrow{OQ}:=\lambda \overrightarrow{OP}+\mu \overrightarrow{OQ}.
\end{array}
\end{equation}

Here we have something similar to the case $m=1$ and  $n=0$.
\begin{lemma}\label{pointsonthesidePQgc}
(i) The coefficient $\lambda$ in (\ref{sideOPnew}) is an integer, and equal to $\widetilde{t}(m-n)$, if

$k=\beta (m^2-mn+n^2)\widetilde{t}+[m\alpha-n(\alpha+\beta)] t$,
$\ell=\widetilde{r} (m^2-mn+n^2)\widetilde{t}-(m\frac{\widetilde{r}+\widetilde{s}}{2}+n\frac{\widetilde{r}-\widetilde{s}}{2})t$, where $t$ and $\widetilde{t}$ are arbitrary integers.
In this case, the coefficient $\mu$ is equal to $n\widetilde{t}-t$.

(ii) The coefficient $\mu$ in (\ref{sideOPnew}) is an integer, and equal to $\widetilde{t}m$, if

$k=-\alpha (m^2-mn+n^2)\widetilde{t}+(m\beta+n\alpha) t$,
$\ell=\frac{\widetilde{r}+\widetilde{s}}{2} (m^2-mn+n^2)\widetilde{t}+(m\widetilde{r}-n\frac{\widetilde{r}+\widetilde{s}}{2})t$, where $t$ and $\widetilde{t}$ are arbitrary integers.
In this case, the coefficient $\lambda$ is equal to $t-n\widetilde{t}$.

(iii) The value  $\mu+\lambda$ is an integer, and equal to $\widetilde{t}m$, if

$k=\beta (m^2-mn+n^2)\widetilde{t}+[m(\alpha+\beta)-n\beta] t$,
$\ell=\widetilde{r} (m^2-mn+n^2)\widetilde{t}-(m\frac{\widetilde{s}-\widetilde{r}}{2}+n\widetilde{r})t$, where $t$ and $\widetilde{t}$ are arbitrary integers.
In this case, the coefficient $\mu$ is equal to $n\widetilde{t}-t$.

(iv) If $\gcd[m\frac{\widetilde{r}+\widetilde{s}}{2}+n\frac{\widetilde{r}-\widetilde{s}}{2},m\alpha-n(\alpha+\beta)]=1$, then there are no points of the lattice in the interior of $\overline{OQ}$.

(v) If $\gcd(m\widetilde{r}-n\frac{\widetilde{r}+\widetilde{s}}{2},m\beta+n\alpha)=1$, then there are no points of the lattice in the interior of $\overline{OP}$.

(vi) If $\gcd(m\frac{\widetilde{s}-\widetilde{r}}{2}+n\widetilde{r} ,m(\alpha+\beta)-n\beta)=1$, then there are no points of the lattice in the interior of $\overline{PQ}$.
\end{lemma}

The proof of this lemma is similar to the proof of Lemma~\ref{pointsonTheSides}. As before, the hypothesis in cases (iv), (v), and (vi) in
Lemma~\ref {pointsonthesidePQgc} can be relaxed.

\begin{lemma}\label{pointsonTheSidesgeneralcase}
(i) If $\gcd(m\frac{\widetilde{r}+\widetilde{s}}{2}+n\frac{\widetilde{r}-\widetilde{s}}{2},m\alpha-n(\alpha+\beta))=\nu$, then the coefficient $\lambda$ is still an integer if $k=\beta (m^2-mn+n^2)\widetilde{t}+(m\alpha-n(\alpha+\beta)) \frac{t}{\nu} $ and
$\ell=\widetilde{r} (m^2-mn+n^2)\widetilde{t}-(m\frac{\widetilde{r}+\widetilde{s}}{2}+n\frac{\widetilde{r}-\widetilde{s}}{2})\frac{t}{\nu}$, where $t$ and $\widetilde{t}$ are arbitrary integers. In this case, the coefficient $\mu$ is equal to $n\widetilde{t}-\frac{t}{\nu}$, and this gives $\nu-1$ points of the lattice in the interior of $\overline{OQ}$.

(ii) If $\gcd(m\widetilde{r}-n\frac{\widetilde{r}+\widetilde{s}}{2},m\beta+n\alpha)=\nu $, then
the coefficient $\mu$ in (\ref{sideOPnew}) is an integer equal to $\widetilde{t}m$, if $k=-\alpha (m^2-mn+n^2)\widetilde{t}+(m\beta+n\alpha) \frac{t}{\nu}$,
$\ell=\frac{\widetilde{r}+\widetilde{s}}{2} (m^2-mn+n^2)\widetilde{t}+(m\widetilde{r}-n\frac{\widetilde{r}+\widetilde{s}}{2})\frac{t}{\nu}$, where $t$ and $\widetilde{t}$ are arbitrary integers.
The coefficient $\lambda$ is equal to $\frac{t}{\nu}-n\widetilde{t}$.

(iii) If $\gcd(m\frac{\widetilde{s}-\widetilde{r}}{2}+n\widetilde{r} ,m(\alpha+\beta)-n\beta)=\nu$, then the value  of $\mu+\lambda$ is  an integer, and equal to $\widetilde{t}m$, if
$k=\beta (m^2-mn+n^2)\widetilde{t}+(m(\alpha+\beta)-n\beta) \frac{t}{\nu}$,
$\ell=\widetilde{r} (m^2-mn+n^2)\widetilde{t}-(m\frac{\widetilde{s}-\widetilde{r}}{2}+n\widetilde{r})\frac{t}{\nu}$, where $t$ and $\widetilde{t}$ are arbitrary integers.
In this case, the coefficient $\mu$ is equal to $n\widetilde{t}-\frac{t}{\nu}$.

(iv) The value $\nu$ in each of the above cases is always a divisor of $d$.
\end{lemma}

\proof. The calculations we have done in showing Lemma~\ref{pointsonthesidePQgc} are valid here with the substitution $\frac{t}{\nu}$ instead of $t$.
For the part (iv), we note that since $\gcd(m,n)=1$ we have $mx+ny=1$ for some
integers $x$ and $y$. Let us consider on (i). One can check the identities
$$\alpha\left(m\frac{\widetilde{r}+\widetilde{s}}{2}+n\frac{\widetilde{r}-\widetilde{s}}{2}\right)-\left(m\alpha-n(\alpha+\beta)\right)\frac{\widetilde{r}+\widetilde{s}}{2}=nd,\ and $$

$$(\alpha+\beta)\left(m\frac{\widetilde{r}+\widetilde{s}}{2}+n\frac{\widetilde{r}-\widetilde{s}}{2}\right)+\left(m\alpha-n(\alpha+\beta)\right)\frac{\widetilde{r}-\widetilde{s}}{2}=md.$$
These two equalities can be combined to get

$$((x+y)\alpha+x\beta)\left(m\frac{\widetilde{r}+\widetilde{s}}{2}+n\frac{\widetilde{r}-\widetilde{s}}{2}\right)+\left(m\alpha-n(\alpha+\beta)\right)(x\frac{\widetilde{r}-\widetilde{s}}{2}-y\frac{\widetilde{r}+\widetilde{s}}{2})=d,$$

\n which implies the claim that $\nu$ must divide $d$.\eproof

\n The following theorem allows one to compute the Ehrhart polynomial for an equilateral triangle in $\mathbb Z^3$.

\begin{theorem}\label{coefficientc1gc}  The Ehrhart polynomial of an equilateral triangle ${{\cal T}^{m,n}_{a,b,c}}$ (given by Theorem~\ref{oldnewparametrization}) with $\gcd(m,n)=1$, is given by
$$L({\cal T}^{m,n}_{a,b,c},t)=\frac{d(m^2-mn+n^2)t^2+c_1t}{2}+1, \ t\in \mathbb N,\ where $$
$$ c_1=\gcd(m\frac{\widetilde{r}+\widetilde{s}}{2}+n\frac{\widetilde{r}-\widetilde{s}}{2},m\alpha-n(\alpha+\beta))+\gcd(m\widetilde{r}-n\frac{\widetilde{r}+\widetilde{s}}{2},m\beta+n\alpha)+\gcd(m\frac{\widetilde{s}-\widetilde{r}}{2}+n\widetilde{r} ,m(\alpha+\beta)-n\beta),$$

\n and ($\alpha$, $\beta$) is a particular solution of the Diophantine equation $\frac{\widetilde{r}+\widetilde{s}}{2}\beta+\widetilde{r}\alpha=d$.
\end{theorem}

\section{The case of $a=b$ and some further questions}

If the  equation of the plane (\ref{planelattice}) has the property that $a=b$ ($2a^2+c^2=3d^2$), then the parametrization of the equilateral triangles given by
Theorem~\ref{oldnewparametrization} simplifies to:

${{\cal T}^{m,n}_{a,b,c}}:=\triangle OPQ$ with
$P$, $Q$ in ${\cal P}_{a,b,c}$, is equilateral if and only if for
some integers $m$, $n$

\begin{equation}\label{vectoridaeqb}
\overrightarrow{OP}=m\overrightarrow{\zeta}-n\overrightarrow{\eta},\
\
\overrightarrow{OQ}=n\overrightarrow{\zeta}+(m-n)\overrightarrow{\eta},
\ { \rm with} \ \overrightarrow{\zeta}=(\zeta_1,\zeta_2,\zeta_3),
\overrightarrow{\eta}=(\eta_1,\eta_2,\eta_3),
\end{equation}

\begin{equation}\label{paramtwoaeqb}
\begin{array}{l}
\begin{cases}
\ds \zeta_1=-\frac{d+c}{2} \\ \\
\ds \zeta_2=\frac{d-c}{2}\\ \\
\ds \zeta_3=a
\end{cases}
 ,\  and \ \ \begin{cases}
\ds \eta_1=\frac{d-c}{2} \\ \\
\ds \eta_2=-\frac{d+c}{2}\\ \\
\ds \eta_3=a
\end{cases}
\end{array}.
\end{equation}

We observe that $\widetilde{r}=\widetilde{s}=1$ and so we can choose $\alpha=d$ and $\beta=0$ to satisfy (\ref{alphabeta}).
Assuming as before that $\gcd(m,n)=1$,  then the Ehrhart polynomial reduces to the simple formula

\begin{equation}\label{aeqb}
L({\cal T}^{m,n}_{a,a,c},t)=\frac{d(m^2-mn+n^2)t^2+[\gcd(m,d)+\gcd(n,d)+gcd(m-n,d)]t}{2}+1, \ t\in \mathbb N.
\end{equation}

In  \cite{ejicrt} we have characterized the primitive triples $(a,c,d)\in \mathbb N^3$ satisfying $2a^2+c^2=3d^2$.
This was done in a manner similar to the way that Pythagorean triples are usually described with a one-to-one correspondence to a special set of pairs of natural numbers.

\begin{theorem}\label{characterizationofx=y}
Suppose that  $k$ and $\ell$ are positive integers with $k$ odd and $\gcd(k,\ell)=1$. Then $a$, $c$ and $d$ given by

\begin{equation}\label{parx=y}
d=2\ell^2+k^2 \ with \ \begin{cases} a=|2\ell^2+ 2k\ell-k^2|, \ c=|k^2+
4k\ell-2\ell^2|,\ if \ k\not\equiv \ell \ {\rm
(mod\ 3)}\\ \\
a=|2\ell^2-2k\ell-k^2|, \ c=|k^2-4k\ell-2\ell^2|,\ if \ k\not\equiv -\ell\  {\rm
(mod\ 3)}
\end{cases}
\end{equation}

\n constitute a positive primitive solution of $2a^2+c^2=3d^2$.

Conversely, with the exception of the trivial solution $a=c=d=1$,
every positive primitive solution for $2a^2+c^2=3d^2$ appears in
the way described above for some $l$ and $k$.
\end{theorem}

In particular, if $d>3$ is a prime of the form $8m+1$ or $8m-5$ ($m\in \mathbb N$), we can find (see \cite{cox})
$k$ and $\ell$ as in Theorem~\ref{characterizationofx=y}. So, we can construct examples as in
the Introduction in which the polynomial is the same for essentially different triangles. As a curiosity, for $d=2011$ we get two solutions
 for $(a,c)$, $(913,3235)$ and $(2461,139)$. The corresponding Ehrhart polynomials are identical:
 $L({\cal T}^{1,0}_{a,a,c},t)=\frac{2011t^2+2013t}{2}+1$, $t\in \mathbb N$.

Our ultimate goal is to compute the Ehrhart polynomials for all of the regular
lattice polyhedrons (tetrahedrons, cubes and octahedrons, see \cite{ejips}), and naturally the first step is to look at their faces.

We pose the following question. What is the cardinality of the set

$$E(d):=\{L({\cal T}^{1,0}_{a,b,c},t) | a^2+b^2+c^2=3d^2\}?$$

This sequence begins in the way recorded in the next table.

\centerline{Table 1}
\vspace{0.1in}
\n
\centerline{\begin{tabular}{|c|c|c|c|}
  \hline   \hline
 d  & All \ primitive\ solutions\ of\ $[a,b,c]$ &  $E(d)$ & $c_1$  \\
   \hline   \hline
  1 & $\{ [1, 1, 1]\}$& 1& $\{3\}$\\ \hline
  3& $\{[1, 1, 5]\}$& 1& $\{5\}$\\ \hline
  5& $\{[1, 5, 7]\}$& 1& $\{3\}$\\
 \hline
 7& $\{[1, 5, 11]\}$& 1& $\{3\}$\\ \hline
9& $\{[1, 11, 11], [5, 7, 13]\}$& 2& $\{5,11\}$\\
\hline
11& $\{[1, 1, 19], [5, 13, 13], [5, 7, 17]\}$& 2& $\{3,13\}$\\
\hline
13& $\{[5, 11, 19], [7, 13, 17]\}$& 1& $\{3\}$\\
\hline
15& $\{[5, 11, 23], [1, 7, 25], [5, 17, 19]\}$& 1& $\{5\}$\\
\hline
17& $\{[7, 17, 23], [1, 5, 29], [13, 13, 23], [11, 11, 25]\}$& 2& $\{3,19\}$\\
\hline
19& $\{[11, 11, 29], [1, 11, 31], [5, 23, 23], [13, 17, 25]\}$& 2& $\{3,21\}$\\
\hline
21& $\{[11, 19, 29], [1, 19, 31], [13, 23, 25]\}$& 1& $\{5\}$\\
\hline
23& $\{[11, 25, 29], [1, 25, 31], [1, 19, 35], [7, 13, 37]\}$& 1& $\{3\}$\\
\hline
25& $\{[1, 5, 43], [17, 25, 31], [11, 23, 35], [5, 13, 41], [17, 19, 35]\}$& 1& $\{3\}$\\
\hline
27& $\{[1, 31, 35], [11, 29, 35], [17, 23, 37], [7, 17, 43], [13, 13, 43]\}$& 3& $\{5,11,29\}$\\
\hline
29& $\{[23, 25, 37], [1, 11, 49], [1, 29, 41], [7, 25, 43], [5, 17, 47]\}$& 1& $\{3\}$\\
\hline
31& $\{[7, 25, 47], [11, 19, 49], [5, 7, 53], [17, 35, 37], [19, 29, 41]\}$& 1& $\{3\}$\\
\hline
33& $\{[5, 29, 49], [7, 37, 43], [23, 37, 37], [19, 35, 41], [25, 31, 41], [23, 23, 47], [13, 17, 53]\}$& 3& $\{5,15,35\}$\\
\hline
35& $\{[5, 29, 53], [17, 19, 55], [11, 23, 55], [25, 37, 41], [25, 29, 47], [5, 13, 59]\}$& 1& $\{3\}$\\
\hline
37& $\{[7, 43, 47], [5, 19, 61], [5, 41, 49], [1, 25, 59], [11, 31, 55], [23, 37, 47]\}$& 1& $\{3\}$\\
\hline
39& $\{[13, 37, 55], [1, 29, 61], [5, 7, 67], [7, 17, 65], [11, 31, 59], [23, 35, 53]\}$& 1& $\{5\}$\\
\hline
41& $\{[5, 47, 53], [1, 1, 71], [19, 31, 61], [13, 43, 55], [25, 47, 47], [5, 23, 67], [31, 41, 49], [17, 23, 65]\}$& 2& $\{3,43\}$\\
\hline \hline
\end{tabular}}

\end{document}